\documentclass[12pt,emtex]{article}
\usepackage{times}
\usepackage{graphicx}
\usepackage{amsfonts} %
\usepackage{amsmath} %[leqno]=Gleichungen werden links numeriert;
\usepackage{amssymb} %
\usepackage{amsthm}%theorem-umgebungen
\usepackage{epsfig}
\usepackage{latexsym}
\usepackage[matrix,curve]{xypic}
\usepackage[german]{babel} %Deutsch
\newcommand{\C}{\Bbb{C}}

\newcommand{\R}{\Bbb{R}}
\newcommand{\Z}{\Bbb{Z}}
\newcommand{\Q}{\Bbb{Q}}
\newcommand{\A}{\Bbb{A}}

\begin{document}
\centerline{\bf The Trace of Hecke operators on the space of classical holomorphic}
\centerline{\bf Siegel modular forms of genus two}

\bigskip\noindent
\centerline{ (Rainer Weissauer) }

\bigskip\noindent
In this note we specialize the results on the trace formula
from [W1],[W2] and [W3] to the case of holomorphic Siegel modular forms of genus
two with special emphasis on the classical  case of forms for
the full Siegel modular group.

\bigskip\noindent
{\it Notation}. Let $(V_\rho,\rho)$ be
an irreducible representation of the linear group
$Gl(g,\C)$. Let $\Gamma$  be a subgroup of finite index of the Siegel modular group
$\Gamma_g = Sp(2g,\Z)$. Let $H_g$ be the Siegel upper half space
of genus $g$, i.e. the space of all complex symmetric $g\times g$ matrices $Z$ with positive definite imaginary part.
Vector valued holomorphic Siegel modular forms
of genus $g$ and type $\rho$ are holomorphic functions $$ f: H_g \to  V_\rho $$
with the transformation property 
$$  f \bigl((AZ+B)(CZ+D)^{-1}\bigr) = \rho(AZ+B) f(Z) $$
for all matrices $$\begin{pmatrix} A & B \cr C & D \cr \end{pmatrix} \ \in \ \Gamma \ . $$
Such a function $f$ is called a cusp form if it is of rapid decay at infinity.
Usually this is expressed in terms of the Fourier expansions of $f$ at the cusps (see [F]).

\bigskip
We restrict now to the case where the genus $g$ is two.
Then we may assume that  $\rho$ is of the form
$\rho(M) = Sym^r(M) \otimes det(M)^k$ for the $r$-th symmetric power
$Sym^r$ of the standard representation
 of $Gl(2,\C)$. In the following let us  assume 
 $k\geq 3$ and we are only interested in the space $[\Gamma,\rho]_0$ of cusp forms within the space $[\Gamma,\rho]$ of all modular forms.
 For $(k_1,k_2)=(r+k,k)$ then $k_1\geq k_2\geq 3$. Then it is well known that a Siegel cusp form $f$ with these properties gives rise to cuspidal automorphic representations
of the adele group $G(\A)=GSp(4,\A)$. Decomposing these representations into a direct sum of irreducible automorphic representations $\Pi = \otimes'_v \Pi_v$
all of the archimedean representations $\Pi_\infty$ which arise 
 belong to the holomorphic discrete series of weight $(k_1,k_2)$
in the sense of [W1].

\bigskip Now we review results of [W1] relevant for our applications.
Let $\A= \R \times \A_{fin}$ denote the ring of rational adeles.  Let $dg$ denote a Haar measure on $G(\A)=GSp(4,\A)$.
The Hilbert space $L^2_0(G(\Q)\setminus G(\A),dg) \subset L^2(G(\Q)\setminus G(\A),dg) $ of cuspidal automorphic representations of $G(\A)$ decomposes discretely into a Hilbert direct sum of irreducible cuspidal automorphic representations $\Pi$ of $G(\A)$.
The space $L^2_0(G(\Q)\setminus G(\A),dg)$ contains the subspaces of CAP-representations $L^2_{CAP}(G(\Q)\setminus G(\A),dg)$ and the subspace $L^2_{endo}(G(\Q)\setminus G(\A),dg)$ of weak endoscopic
lifts. The intersection of these two subspaces is zero. See [W1], page 70 resp. [S]
for further details. For the following it suffices to know, that CAP representations 
are the irreducible cuspidal representations, which are weakly equivalent to constituents of
globally induced automorphic representations. Notice that two irreducible  automorphic representations $\Pi_1,\Pi_2$ of $GSp(4,\A)$ are said to be weakly equivalent if their local components $\Pi_{1,v},\Pi_{2,v}$ are isomorphic for almost all places $v$.
This being said let $L^2_{00}(G(\Q)\setminus G(\A),dg)$ denote the orthogonal
complement of these two subspaces in $L^2_0(G(\Q)\setminus G(\A),dg)$.
An irreducible constituent $\Pi=\Pi_\infty \otimes \Pi_{fin}$ of $L^2_{0}(G(\Q)\setminus G(\A),dg)$
is said to be a cohomological representation, if its archimedean component
$\Pi_\infty$ belongs to the discrete series representations of the group $GSp(4,\R)$. This condition is equivalent to the condition that there exist integers $k_1\geq k_2 \geq 3$
such that $\Pi_\infty $ belongs to a local archimedean $L$-packet $\{ \Pi_\infty^{hol}, \Pi_\infty^{W} \}$ of cardinality  two attached
to this weight $(k_1,k_2)$. 
For cohomological irreducible cuspidal automorphic representations $\Pi$ not of CAP-type
we constructed in [W1] associated four dimensional Galois representations $\rho_{\Pi,\lambda}$ of the absolute Galois group
of $\Q$ with coefficients in the algebraic closure $\overline{\Q_l}$ of $\Q_l$, 
which are defined over some finite dimensional extension field $E_\lambda$ of the $l$-adic field $\Q_l$.
Once and for all we fix a field isomorphism $ \tau: \overline{\Q_l} \cong \C$
and  tacitly identify $\overline{\Q_l}$ with the field of complex numbers.

\bigskip
The Shimura variety $M=GSp(4,\Q)\setminus GSp(4,\A)/ K_\infty$ has a model over the reflex field $E=\Q$. Here $K_\infty \subset GSp(4,\R)$ denotes the stabilizer of the point $i\cdot E$ in $H_2$ so that $GSp(4,\R)/K_\infty$ can be identified with the union of half spaces $H_2 \cup - H_2$.
To a representation $\rho$ of $Gl(2,\C)$ as above one can attach a $\overline{\Q_l}$-coefficient system $V_\lambda$ for $\lambda=\lambda(\rho)$ on $M$
and decompose the etale cohomology $H^\bullet_c(M, V_\lambda)$ 
as a representation of the group $ GSp(4,\A_{fin})  $. It is known that this representation is automorphic. On the Eisenstein cohomology
$$ H_{Eis}^\bullet(M,V_\lambda) = Kern(H^\bullet_c(M,V_\lambda) \to H^\bullet(M,V_\lambda)) $$
the group $GSp(4,\A_{fin})$ acts with constituents of globally induced representations.
The image $H_!^\bullet(M,V_\lambda) = Im(H^\bullet_c(M,V_\lambda) \to H^\bullet(M,V_\lambda)) $
is completely decomposable into irreducible representations of $GSp(4,\A_{fin})$ ([H]).
$H_!^\bullet(M,V_\lambda)$ contains a maximal subspace $H^\bullet_E(M,V_\lambda)$, whose irreducible constituents are weakly equivalent to 
globally induced representations. Using transcendent methods one can show $H^\bullet_E(M,V_\lambda) = H^\bullet_{res}(M,V_\lambda) \oplus H^\bullet_{CAP}(M,V_\lambda)$. Classes in the first subspace are represented by residues of Eisenstein series. The second subspace is the part of the cuspidal cohomology defined by the  CAP-representations. The orthocomplement of $H^\bullet_E(M,V_\lambda)$ in $H^\bullet_!(M,V_\lambda)$
with respect to the cup-product decomposes discretely as a module under the group $Gal(\overline{\Q}:\Q) \times GSp(4,\A_{fin})  $ 
$$   \bigoplus_{\Pi} \tilde\rho_{\Pi,\lambda} \otimes  \Pi_{fin}  \ .$$
where the summation extends over all irreducible cohomological cuspidal automorphic representations $\Pi = \Pi_\infty \otimes \Pi_{fin}$ of $G(\A)$ not of CAP-type. 
It can again be  split up into two subspaces.
One is the  subspace $H^\bullet_{endo}(M, V_\lambda)$ defined by the  weak endoscopic
lifts $\Pi$, the other is  the subspace $H^\bullet_{00}(M, V_\lambda)$ defined by the representations $\Pi$  in $L^2_{00}(G(\Q)\setminus G(\A),dg)$. The nature of the Galois representation $\tilde\rho_{\Pi,\lambda}$ 
depends on the type of $\Pi$ in this sense.
The \lq{motivic}\rq\ Galois representations
$\tilde\rho_{\Pi,\lambda}$ of $Gal(\overline{\Q}:\Q)$  are of finite dimension over $\overline{\Q_l}$ and they are uniquely
determined by the weak equivalence class of the automorphic representation $\Pi$. This  easily follows from  
the Cebotarev density theorem. 
 
 \bigskip\noindent
We will not discuss the Eisenstein cohomology, which is explained in greater detail in [H],[P],[FG] and [BFG]. Our main focus will be on the cases where $\Pi$ is a weak endoscopic lift or belongs to
$L^2_{00}(G(\Q)\setminus G(\A),dg)$. For the latter case we have the following

\bigskip\noindent
{\bf Theorem 1} (Stability). {\it Suppose $\Pi=\Pi_\infty\otimes \Pi_{fin}$ is an irreducible cuspidal automorphic representation in $L^2_{00}(G(\Q)\setminus G(\A),dg)$ for which
$\Pi_\infty$ is in the local archimedean  $L$-packet  $\{ \Pi_\infty^{hol}, \Pi_\infty^{W} \}$ of a discrete
series representation of weight $(k_1,k_2)$. Then the multiplicities of the representations $\Pi_{fin}\otimes  \Pi_\infty^{hol}$ and $\Pi_{fin}\otimes \Pi_\infty^{W}$ in $H^\bullet_{00}(M, V_\lambda)$, or equivalently multipicities in $L^2_{00}(G(\Q)\setminus G(\A),dg)$, coincide
$$  m(\Pi_{fin}\otimes  \Pi_\infty^{hol} )\ = \  m(\Pi_{fin}\otimes \Pi_\infty^{W}) \ .$$
The semisimplification $ \tilde\rho_{\Pi,\lambda}^{ss}$ of the motivic representation $\tilde\rho_{\Pi,\lambda}$ is concentrated in the cohomology of  degree three and 
 and is isomorphic to an isotypic
multiple $$ \tilde\rho_{\Pi,\lambda}^{ss}=  n(\Pi)  \cdot \rho_{\Pi,\lambda} $$ of the four-dimensional symplectic Galois representations $\rho_{\Pi,\lambda}$ attached to $\Pi$. It is defined over
a finite extension field $E_\lambda$ of $\Q_l$. Viewed as a representation over $\Q_l$ the representation  $\rho_{\Pi,\lambda}$ is of Hodge-Tate type, and  
its Hodge-Tate components $(k_1+k_2-3,0),(k_1-1,k_2-2),(k_2-2,k_1-1),(0,k_1+k_2-3)$
occur with the same multiplicity $[E_\lambda:\Q_l]$.}

\bigskip\noindent
{\it Proof}. By [W3], theorem 1 any irreducible cuspidal automorphic  $\Pi$ of $GSp(4,\A)$, with the assumptions as in the theorem above, is weakly equivalent to a globally generic cuspidal automorphic representation $\Pi'$ of $GSp(4,\A)$ for which the local archimedean component $\Pi'_\infty$ is in the same local archimedean local $L$-packet as $\Pi_{\infty}$. This assertion allows to apply [W1], theorem III and [W1], proposition 1.5., which immediately give the statements of the theorem above.  It should be remarked that the results proven in [W1] depend on certain hypotheses A and B
made in loc. cit on page 70 and page 80. The proof of the hypotheses A and B is 
the main content of [W2]. QED

\bigskip
Now we apply the last theorem.
The case of our particular interest is the case of the full Siegel modular group $$\Gamma = \Gamma_{2} \ .$$  Let $[\Gamma_2,\rho]_0$ be the corresponding space of vector valued
holomorphic Siegel cusp forms for which the weight $(k_1,k_2)$ of $\rho$ satisfies $k_1\geq k_2\geq 3$. Under the action of the algebra ${\cal H}$ of spherical Hecke operators $T\in {\cal H}$ (see [F]) every  cuspform $f$ in $[\Gamma_2,\rho]_0$ can be decomposed into a finite sum of eigenforms of ${\cal H}$. For a cusp form $f$, which is an eigenform of all Hecke operators $T$
$$   f\vert_\rho T = \lambda(T) \cdot f \ ,$$
the eigenvalues $\lambda(T)\in \C$ define an algebra homomorphism $$ \lambda_f : {\cal H} \to \C \ .$$

\bigskip\noindent
{\bf Theorem 2} (Multiplicity one). {\it For the case of the full Siegel modular group $\Gamma_2$ the homomorphism $\lambda_f$ uniquely determines the weight $k_1,k_2$ of
$\rho$ and uniquely determines the eigenform $f\in [\Gamma_2,\rho]_0$ up to a scalar.}   

\bigskip\noindent
{\it Proof}. Any cuspidal eigenform $f$ of ${\cal H}$ determines
an irreducible cuspidal automorphic representation $\Pi=\Pi(f)$ for which
$\Pi_\infty(f)$ belongs to the holomorphic discrete series $\Pi_\infty^{hol}$ of type $(k_1,k_2)$. Conversely
any irreducible cuspidal automorphic representation with $\Pi_\infty$ in the holomorphic discrete series of type $(k_1,k_2)$ determines a holomorphic cuspidal Siegel eigenform of all Hecke
operators by considering the one dimensional space of  spherical vectors in $\Pi_{fin}$ and the one dimensional space of lowest $K_\infty$-type in $\Pi_\infty$.   Notice that $\lambda_f$ determines the degree four $L$-series of $f$ or $L(\Pi_{fin},s)=L(f,s)$,
and also the degree five $L$-series $\zeta(\Pi,s)$.
Conversely $L(f,s)=L(\Pi_{fin},s)$ determines the spherical representation $\Pi_{fin}$.
Therefore, 
since the CAP-cases are characterized by poles of their degree four or degree five $L$-series (see [PS],[S]), the CAP-property is detected by $\lambda_f$. So this CAP-case can be dealt with separately. In fact in the CAP-case the statement reduces to a statement on forms in the Maass Spezialschar, where this is well known ([PS], [Z]).
So we may assume without restriction of generality that either $f$ defines a weak endoscopic lift $\Pi(f)$ or a representation $\Pi(f)$ of $L^2_{00}$-type. In both these cases $\Pi (f)$ is weakly equivalent to a globally generic representation $\Pi'(f)$ whose archimedean component
$\Pi'_\infty(f)$ belongs to the same local archimedean $L$-packet of weight $(k_1,k_2)$. For the $L^2_{00}$-case this follows from [W3] as already explained. In the case of a weak endoscopic lift $\Pi(f)$ in $\Pi(\sigma)$ this is shown in [W2] Theorem 5.2, page 186. In fact the multiplicity formula of loc. cit. theorem 5.2.4 implies
 $m(\Pi')=1$ for $\Pi'=\otimes'_v \Pi_+(\sigma_v)$. The representation $\Pi'$ is weakly equivalent to $\Pi$ and it is globally generic (see [W2] theorem 4.1 and 4.2 and the references given there). The detailed description of the local representations $\Pi_+(\sigma_v),\Pi_-(\sigma_v)$
given in [W2]  moreover implies  for the full Siegel modular group  that
weak endoscopic lifts do not occur in the space of holomorphic vector valued cusp forms of weight $k_1 \geq k_2\geq 3$.  
This is  discussed in lemma 1 below. Using this we may therefore assume that we are in the $L^2_{00}$-case. This allows us to apply our theorem 1:

\bigskip\noindent
{\it Step 1}. $\lambda_f$ determines $k_1,k_2$. As explained above
$\lambda_f$ determines the partial $L$-series $$L(f,s)=L(\rho_{\Pi,\lambda},s)$$
of  the automorphic representation $\Pi=\Pi(f)$. Hence $\lambda_f$ determines the Galois representation $\rho_{\Pi,\lambda}$ attached to $f$ by the Cebotarev density theorem. Since this representation $\rho_{\Pi,\lambda}$, considered as a representation over $\Q_l$ is a Hodge-Tate Galois representation, we can consider its Hodge-Tate decomposition. The  Hodge-Tate decomposition, described in theorem 1, obviously determines the integers
$k_1$ and $k_2$. 

\bigskip\noindent
{\it Step 2}. The multiplicity one statement. To show it we apply [W1], lemma 1.2
for the two weakly equivalent representations $\Pi = \Pi(f)$ and the globally
generic representation $\Pi'$ associated to it. We use that $\Pi_{fin}=\Pi_{fin}(f)$ is a spherical
representation, since $f$ is a cusp form for the full Siegel modular group.
Hence [W1], lemma 1.2 implies, that also $\Pi'_{fin}$ has to be spherical and moreover that
$$   \Pi_{fin}(f) = \Pi'_{fin}(f)  $$
holds. Let me briefly remark, that this argument uses the global functional equation of the $L$-series attached to 
cuspidal automorphic representations of $GSp(4,\A)$.
 On the other hand
$$ \Pi_\infty (f)= \Pi^{hol}_\infty \quad , \quad \Pi'_\infty(f)= \Pi^{W}_\infty \ , $$
since $f$ is holomorphic and since the generic representation $\Pi'_\infty$ 
has a Whittaker model.
Hence $\Pi_\infty(f) \neq \Pi'_\infty(f)$, but they are contained in the same local archimedean local
$L$-packet attached to $(k_1,k_2)$. 
So we can apply the stability theorem 1 from above, since we are in case $L^2_{00}$. This gives the following multiplicity formula
$$ m(\Pi(f)) = m(\Pi_\infty^{hol}\otimes \Pi_{fin}) = m(\Pi_\infty^{W}\otimes \Pi_{fin}) = m(\Pi') \ .$$ 
But globally generic automorphic representations $\Pi'$  have multiplicity
$m(\Pi')=1$ in the cuspidal spectrum as shown in [JS]. This proves $m(\Pi(f))=m(\Pi')=1$
and gives the second assertion of theorem 2. QED.

\bigskip
For the proof of theorem 2 we still have to show 

\bigskip\noindent
{\bf Lemma 1}. {\it For the full Siegel modular group $\Gamma_2$ the subspace generated by cuspidal eigenforms $f$ of $[\Gamma_2,\rho]_0$, for which $\Pi(f)$ is a weak endoscopic lift, is zero.}

\bigskip\noindent
Recall that by definition  (see [W1] page 70)
a global weak endoscopic cuspidal lift $\Pi=\Pi(\sigma)$ is attached to a pair of holomorphic cuspidal elliptic eigenforms $(f_1,f_2)$ respectively the pair of automorphic irreducible cuspidal
representations $\sigma=(\pi_1,\pi_2)$ of $M(\A)=Gl(2,\A)\times Gl(2,\A)$ associated to the forms $f_1,f_2$ (having the same nebentype character). Let us fix $(f_1,f_2)$ or equivalently  $\sigma=(\pi_1,\pi_2)$.
If $f$ is a weak endoscopic lift, then for a finite set $S$ of exceptional places $$ L^S(f,s)= L^S(f_1,s)L^S(f_2,s)$$ holds for the partial $L$-series $L(f_i,s)$ of the two elliptic cusp forms $f_1,f_2$.
This uniquely determines $\Pi_v=\Pi_v(f)$ outside a finite set  of places $S$.
We also know $\Pi_\infty(f) = \Pi_\infty^{hol}$ for the holomorphic Siegel cusp form $f$.
For a cuspidal weak endoscopic lift $\Pi$ the local components
$\Pi_v$ at the nonarchimedean places have been described in [W2]: 
For the places $v\in S,v\neq \infty$ 
either $\Pi_v$ is uniquely determined
$$  \Pi_v \in \{ \Pi_+(\sigma_v) \} $$
(i.e. the local $L$-packet of the lift attached to $\sigma$ has cardinality one) or alternatively 
there are two possible choices in the local $L$-packet determined by the lift
of the local representation $\sigma_v=(\pi_{1,v},\pi_{2,v})$ 
$$ \Pi_v \in \{ \Pi_+(\sigma_v), \Pi_-(\sigma_v) \} \ .$$

\bigskip\noindent
{\it Example} (see [W2] page 153). $\Pi_+(\sigma_\infty)=\Pi_\infty^W$ and $\Pi_-(\sigma_\infty)=\Pi_\infty^{hol}$
in the archimedean local $L$-packet of $GSp(4,\R)$ defined by the discrete series
representation $\sigma_\infty$ of $M(\R)$.

\bigskip\noindent
For a more detailed description of these local $L$-packets and the proofs we refer to [W2], section 4.11.  See loc. cit. page 153 for an overview,  and also [W2], theorem 5.2.  A brief review of the main results can be found in the formulation of hypotheses A in [W1]. 
Unfortunately  the formulation of hypotheses A, part (4) is misstated. It should read: \lq{For the finitely many places $v$ of $F$ for which $\sigma_v$ belongs to the discrete series
of the group $M(F_v)$ the representation $\Pi_v$ is contained in a local $L$-packet  $\{\Pi_+(\sigma_v), \Pi_-(\sigma_v) \}$ consisting of two classes of irreducible admissible representations $\Pi_{\pm}(\sigma_v)$ of $GSp(F_v)$, which only depend on $\sigma_v=(\pi_{1,v},\pi_{2,v})$. At the
remaining places, where $\sigma_v$ does not belong to the discrete series, $\Pi_v \cong \Pi_+(\sigma_v)$ is uniquely determined by $\sigma_v$}\rq.  This being said let us describe
 the global picture. The main global result  is the following.
Let $\Sigma$ be the set of places for which the local component $\sigma_v$ of $\sigma=\otimes'_v \sigma_v$ is in the discrete series.
For $v \in \Sigma$ fix signs $\varepsilon_v = \pm 1$. 
By a slight abuse of notation we now write $\Pi_{\varepsilon_v}(\sigma_v)$
with $\varepsilon_v \in \{\pm 1\}$ instead of using the indices $ \pm $. Then, with this convention,
the irreducible representation
$$ \bigotimes_{v\in\Sigma} \ \Pi_{\varepsilon_v}(\sigma_v) \ \otimes \ \bigotimes_{v\notin \Sigma}{}'
\ \Pi_{+}(\sigma_v) $$
appears with the multiplicity $\frac{1}{2}(1+ \prod_{v\in \Sigma} \varepsilon_v)$ in the cuspidal spectrum ([W2], theorem 5.2.4).
Hence the multiplicity is zero or one depending on whether $\prod_{v\in \Sigma} \varepsilon_v$ is equal to $-1$ or $1$.  In our case this gives 
 
\bigskip\noindent
{\it Proof of lemma 1}. Assume there exists  a holomorphic Siegel cusp form $f$ for the full modular group in the lift $\Pi(\sigma)$ for some global $\sigma$ as above. The associated cuspidal automorphic representation $\Pi=\Pi(f)$ has a  spherical nonarchimedean  representation $\Pi_{fin}(f)$, and it occurs in the lift $\Pi(\sigma)$. Hence $$ \Pi_v(f) = \Pi_{\varepsilon_v}(\sigma_v) \quad , \quad \varepsilon_v \in \{\pm 1 \} $$
for the finitely many places $v\in\Sigma$ where $\sigma_v$ belongs to the discrete series.
Checking the list of the possibilities for  
 $\Pi_v(f) = \Pi_{\varepsilon_v}(\sigma_v)$ for nonarchimedean $v\in \Sigma$ in [W2], page 153
 we see that either $\Pi_v(f)$ has to be in the discrete series or has to be a limit of discrete series
 depending on whether $\sigma_v^*\cong \sigma_v$ or not. 
On the other hand $\Pi_v=\Pi_v(f)$ is spherical for all nonarchimedean places $v$.
Hence $\Pi_v$ can not be a limit of discrete series or a discrete series representation
for any nonarchimedean place $v$. This implies
$$  \Sigma = \Sigma(\sigma) = \{ \infty \} \ .$$ 
In other words only the archimedean component $\sigma_\infty$ belongs to the discrete series.
But then the multiplicity formula above implies that $\varepsilon_\infty = \prod_{v\in \Sigma} \varepsilon_v = 1$, or in other words
$$ \Pi_\infty(f) = \Pi_+(\sigma_\infty) \ .$$
Since $\Pi_+(\sigma_\infty)$ is the representation $\Pi_\infty^W$ with a Whittaker
model this implies $\Pi_\infty(f) =\Pi_\infty^W$ contradicting
the fact that $\Pi_\infty(f)=\Pi^{hol}_{\infty}$ (holomorphicity of $f$).
This contradiction proves the lemma. QED

\bigskip
We remark, that in the proof above one could alternatively use the fact that spherical nonarchimedean representations $\Pi_v$ have a Whittaker model to avoid the argument
with (limits of) discrete series representation for $\Pi_v$.

\bigskip\noindent
{\it Remark}. As a consequence we see that for the case of the full modular group
all the contribution of weak endoscopic lifts $\Pi$ for a fixed cuspidal $\sigma$ 
to the representation to the cohomology groups of the Siegel modular variety
$$  {\cal A}_2 = GSp(4,\Q)\setminus GSp(4,\A) / (K_\infty \times GSp(4,\Z_{fin})) $$
is restricted to the nonholomorphic cohomology of Hodge type
$(k_1-1,k_2-2)$ and $(k_2-2,k_1-1)$  in the cohomology degree three. 
So it remains to discuss the motivic Galois representation attached to this nonholomorphic contribution of the 
lift $\Pi$ in the the weak endoscopic lift of each cuspidal $\sigma$.

\bigskip
In general, for a weak endoscopic lift the $\Pi_{fin}$-isotypic component
$$ \tilde\rho_{\Pi} \otimes \Pi_{fin}  \subset H^\bullet_{endo}(M,V_\lambda) $$
has been computed in [W2]. See assertion (7) on page 71 of [W1]. Of course $V_\lambda$ here is determined by $\sigma_\infty$ and vice versa.

Let us restrict this general statement to the case of the full modular group.
Since we consider the full modular group
the representation $\Pi_{fin}$ has to be spherical.
Then we  know that $\Pi$ gives a contribution to the cohomology
$H^\bullet_{endo}({\cal A}_2,V_\lambda)$ only if the following holds 
$$ (*) \quad  \quad \Pi = \Pi^W_\infty \otimes \Pi_{fin}(\sigma) \quad , \quad \Pi_{fin}(\sigma) \mbox{ is spherical} \ $$
as shown in the proof of lemma 1. 

\bigskip
For $\Pi$ as in (*), or  for weak lifts in general, the cohomological trace formula (see [W2] p.81, [W2] section 4.3)  computes the motivic Galois representation in terms of the two-dimensional representations $\rho_1$ resp. $\rho_2$  attached to the $\pi_1$ resp. $\pi_2$  (see ([D]) with certain multiplicities $m_1,m_2$
$$ (\tilde\rho_\Pi)^{ss} \ = \ m_1 \cdot \rho_1 \oplus m_2 \cdot (\rho_2 \otimes \nu_l^{k_2-2}) $$
 that are computable in terms of Hodge theory ([W2] corollary 4.1 and corollary 4.4). By the formula of [W2] at the bottom of page 88 one has
$$ m_1 = m(\Pi^-(\sigma_{\infty})\otimes \Pi_{fin}) = m(\Pi^{hol}(\sigma_{\infty})\otimes \Pi_{fin})$$
and
$$ m_2 = m(\Pi^+(\sigma_{\infty})\otimes \Pi_{fin}) = m(\Pi^{W}(\sigma_{\infty})\otimes \Pi_{fin}) \ .$$
In the  case (*) relevant for the full Siegel modular  group everything simplifies. Indeed by the  proof of lemma 1 we have already seen that
$$ m_1 = 0 \quad , \quad m_2 = 1 \ .$$
Hence $ (\tilde\rho_\Pi)^{ss} =  \rho_2 \otimes \nu_l^{k_2-2} $
is the two dimensional  $\overline{\Q_l}$-adic representation attached to $f_2$ by Deligne [D]. Since these representations are irreducible, as shown by Ribet [R], we get $ (\tilde\rho_\Pi)^{ss} = \tilde\rho_\Pi$, hence
$ \tilde\rho_\Pi = \rho_2 \otimes \nu_l^{k_2-2} $.
So we obtain
$$ H^\bullet_{endo}(M_K, V_\lambda) = H^3_{endo}(M_K, V_\lambda)\ =\  \bigoplus_{\Pi}  \ (\rho_2 \otimes \nu_l^{k_2-2}) \otimes (\Pi_{fin})^K  $$
where the sum runs over all endoscopic lifts $\Pi$ for all $\sigma$ with $\sigma_\infty$ fixed and determined by $\lambda$ or $(k_1,k_2)$. Notice that $(f_1,f_2)$ are distinguished by their weights $$ r_1=k_1+k_2-2 \ \ \mbox{ and } \ \ \ r_2=k_1-k_2+2 \ .$$ In other words $f_2$ is the form of the pair $(f_1,f_2)$ with the lower weight
$r_2 < r_1$. See [W2] page 64 and 77, and also [W2] page  289 or for a brief overview
[W1], page 70. 

\bigskip
It remains to discuss the relevant representations $\sigma$. For this still assume that we consider the case of the Siegel modular variety $M_K={\cal A}_2$ for the maximal compact group $K=GSp(4,\Z_{fin})$ in $GSp(4,\A_{fin})$.
Then $\Pi_{fin}^K$ is one dimensional if $\Pi_{fin}$ is spherical, and $\Pi_{fin}^K$ is zero
otherwise. This and the multiplicity formula $m_2=1,m_1=0$ from above implies that the representation
$$  \sigma = \sigma_\infty \otimes \sigma_{fin} $$
is uniquely determined by $\Pi_{fin}$ and vice versa.
Indeed $\sigma_\infty$ is determined by the coefficient system $V_\lambda$,
and $\sigma_{fin}$ is determined by the lifting formula [W2], lemma 4.27.
In fact by the proof of lemma 1 we know that $\sigma_v$ can not be in the discrete
series, if a spherical representation $\pi_v$ is in its local weak endoscopic $L$-packet
of $\sigma_v$. Hence $\sigma_v$ must be an induced representation by the local classification theory of admissible irreducible representations. But for an induced representation
$\sigma_v$ the unique local endoscopic lift $\Pi^+(\sigma_v)$ again is induced.
It is described by the formula of [W2] lemma 4.27. This formula moreover implies
that $\sigma_v$ is spherical, if $\Pi_v$ is spherical. Conversely, if $\sigma_v$ is spherical,
then also the unique endoscopic lift $\Pi_v=\Pi^+(\sigma_v)$ is spherical. Since $\sigma_{fin}$ is spherical this describes the possible irreducible representations $\rho$ in terms of pairs of classical elliptic cusp forms $f_1,f_2$ of weight $r_1,r_2$ respectively.

\bigskip\noindent
{\bf Lemma 2}. {\it Fix weights $k_1\geq k_2\geq 3$ and a corresponding
coefficient system $V_\lambda$. Let $r_1=k_1+k_2-2$ and $r_2=k_1-k_2+2$. Then the $\overline{\Q_l}$-adic Galois representation of $Gal(\overline \Q: \Q)$ on the cohomology group
$$ H^\bullet_{endo}({\cal A}_2, V_\lambda)  := H^\bullet_{endo}(M, V_\lambda)^K \quad , \quad K=GSp(4,\Z_{fin}) $$
is nontrivial only in cohomology degree three, where it is isomorphic as a representation 
of $Gal(\overline \Q : \Q)$ to $$ 
  \bigoplus_{\sigma_2} \ \dim_\C([\Gamma_1,r_1]_0) \cdot (\rho_2 \otimes \nu_l^{k_2-2})
 $$
with summation  over all two-dimensional $\overline{\Q_l}$-adic  Galois representations $\rho_2$ attached
to the elliptic cuspidal eigenforms $f_2\in [\Gamma_1,r_2]_0$.}

\bigskip
By our definitions we decomposed the Euler chacteristics
 $$e_c({\cal A}_2,V_\lambda) = \sum_{i=0}^6 (-1)^i [H^i_c({\cal A}_2,V_\lambda)] $$
$$
= e_{Eis}({\cal A}_2,V_\lambda) + e_{E}({\cal A}_2,V_\lambda) + e_{endo}({\cal A}_2,V_\lambda) +
e_{00}({\cal A}_2,V_\lambda)\  $$
where the endoscopic term $ e_{endo}({\cal A}_2,V_\lambda) $ 
is given by lemma 2 (up to an additional sign $(-1)^3$ from the 
Euler characteristics). The term $ e_{00}({\cal A}_2,V_\lambda) $ corresponds up to a sign $(-1)^3$ to a \lq{motif\rq}\ of rank $4\cdot dim_{\C}([\Gamma_2,\rho]_{00})$ by theorem 1.   
Together with the next lemma this confirms conjecture 4.1 of [FG]
(notice $m=k_2-3$ and $l-m=k_1-k_2$  in the notation
of loc. cit. )  

\bigskip\noindent
{\bf Lemma 3}. {\it 
$H^\bullet_{E}(M,V_\lambda)$ vanishes for regular weights $k_2>k_1 >3$.}

\bigskip
In other words we obtain the following
formula for the trace of the Hecke operator $T(p)$
for prime $p$ in terms of traces of the Frobenius $F_p$ on cohomology.

\bigskip\noindent
{\bf Theorem 3}. {\it  For regular weights $k_1>k_2>3$ and the corresponding coefficient system $V_\lambda$ on ${\cal A}_2$
the trace of the Hecke operators $T=T(p)$ on the space of holomorphic Siegel cusp forms $[\Gamma_2,\rho]_0$ for the full Siegel modular group $\Gamma_2$ is given by $$  4\cdot trace(T(p),[\Gamma_2,\rho]_0) =  -trace\bigl(F_p,[H_c^\bullet({\cal A}_2,V_\lambda)]\bigr) $$
 $$ - \sum_{\sigma_2} \ \dim_\C([\Gamma_1,r_1]_0) \cdot trace(F_p,\rho_2 \otimes \nu_l^{k_2-2})+ e_{Eis}({\cal A}_2,V_\lambda)  \ .$$}

\bigskip\noindent
{\it Remark}.  Using the formula for $e_{Eis}({\cal A}_2,V_\lambda)$ in [FG] the terms in the second row of the formula of theorem 3 can be expressed in the form
$$ \ \dim_\C\bigl([\Gamma_1,r_1]_0\bigr) \cdot trace\bigl(F_p,[H_c^\bullet({\cal A}_1,V_\mu)])\otimes \nu_l^{k_2-2}\bigr) \  +  \  \dim_\C([\Gamma_1,r_2]_0) $$
$$ + (-1)^{k_1}  \cdot trace(F_p,[H_c^\bullet({\cal A}_1,V_{\mu'})]))  + (1 + (-1)^{k_1})/2 \ .$$  
Here $V_\mu$ is the $\overline \Q_l$-coefficient system
on ${\cal A}_1$ whose cohomology is related to the elliptic modular forms $[\Gamma_1,r_2]$
of weight $r_2$ (by the Eichler-Shimura isomorphism) respectively $V_{\mu'}$ 
is the $\overline \Q_l$-coefficient system
on ${\cal A}_1$ whose cohomology is related to the elliptic modular forms $[\Gamma_1,k]$
of weight $k=k_1$ (if $k_1$ is even) or $k=k_2-1$ (if $k_1$ is odd).
%A similar formula should still be true for $k_1>k_3=3$. 
%See the comment [FG] \S 2.
%Since there are no CAP-contributions in these cases either, the proposition above extends
%to $k_3=3$.

%\bigskip\noindent
%{\it Remark}. 
%The cohomological $trace(F_p,[H_c^\bullet(M,V_\lambda)])$  has been computed by Kottwitz [K]. See also the formula
%[W2].

\bigskip\noindent
{\it Remarks on level 2}. In [BFG] certain explicit formulas were conjectured
for the principal congruence group $\Gamma=\Gamma_2[2]$ of level 2.
For example the first part of conjecture 7 amounts  
to a certain property of the 2-adic representations $\Pi_v=\Pi_\pm(\sigma_v)$ for the 2-adic field $F_v=\Q_2$, namely 
 that locally at the 2-adic place $v$
$$  dim_\C(\Pi_-(\sigma_v)^K) = 1 \quad \mbox{for} \quad \chi_{1,v}/\chi_{2,v} = \chi_0 $$
$$  dim_\C(\Pi_-(\sigma_v)^K) = 5 \quad \mbox{for} \quad \chi_{1,v}/\chi_{2,v} = 1 $$
holds for the 2-adic principal congruence group $K\subset GSp(4,\Z_2)$ of level two in the case
where $\sigma_v=(\pi_{1,v},\pi_{2,v})$ are special representations
$\pi_{i,v} =Sp \otimes \chi_{i,v}$ of $Gl(2,\Q_2)$ whose character $\chi_{i,v}$ is either trivial or equal to the nontrivial unramified quadratic character $\chi_0$. Implicitly in the regular   case there is even the stronger  conjecture 7.4 of [BFG] that $\Pi_-(\sigma_v)^K$, as module under the symmetric group $GSp(4,\Z_2)/K\cong \Sigma_6$, is isomorphic to $s[1^6]$ resp. $s[2^3]$. To show this conjecture, one has to compute
the local representations $\Pi_-(\sigma_v)$ as in  [W2], case 1c respectively case 1d (page 129f).
As these statements are of local nature, one could prove them locally. However, for a proof it suffices to know that they hold in a single global example.

\bigskip\noindent
{\it Sketch of proof for lemma 3}. Regularity implies the vanishing of cohomology outside
of degree 3 and $H_!^3(M,V_\lambda)=H^3_{cusp}(M,V_\lambda)$ (see [T], p.294).
Hence $H^\bullet_E(M,V_\lambda)= H^\bullet_{CAP}(M,V_\lambda)$ and only representations
contribute with $\Pi_\infty$ in the discrete series. But CAP easily implies regularity by computing the well known archimedean theta lifts defining the CAP representations [S], [PS].

% We use that for endoscopic lifts $\Pi$ in $\Pi(\sigma)$ for $\sigma=(\pi_1,\pi_2)$ the central %character $\omega_\Pi$ and the central characters $\omega_1,\omega_2$
%of $\pi_1,\pi_2$ all coincide. See []. For representations $\Pi$ of level 2 therefore the central %character
%$\omega_\Pi$ is a Dirichlet character of conductor at most 4.  

\bigskip\noindent

\bigskip\noindent

\goodbreak
\bigskip\noindent\centerline{\bf References}

\bigskip\noindent
[D] Deligne P., Formes modulaires et representations l-adiques, In:
Sem. Bourbaki 1968/69, Springer Lecture Notes 79,  pp. 347-363 (1971) 

\bigskip\noindent
[F] Freitag E., Siegelsche Modulfunktionen, Grundlehren der mathematischen
Wissenschaften 254, Springer (1983)

\bigskip\noindent
[FG] Faber C.-van der Geer G., Sur la cohomologie des systemes locaux sur les espaces
des modules des courbes de genre 2 et des surfaces abeliennes, I, II C.R.Acad. Sci. Paris, Ser. I, 338 (2004), 381 - 384, 467 - 470 

\bigskip\noindent
[BFG] Bergstr\"om J.- Faber C.- vand der Geer G., Siegel modular forms of genus 2 and level 2:
cohomological computations and conjectures, International Mathematics Research Notices, vol. 2008.  

\bigskip\noindent
[H] Harder G., Eisensteinkohomologie und die Konstruktion gemischter Motive,
Springer Lecture Notes in Mathematics  1562 (1993)

\bigskip\noindent
[K] Kottwitz R.E., Shimura varieties and $\lambda$-adic representations, In: Ann Arbor Proceedings,
edited by Clozel L., Milne J.S, Perspectives of mathematics, pp 161 - 209 (1990) 

\bigskip\noindent
[P] Pink R., On $l$-adic sheaves on Shimura varieties and their higher
direct images  in the Baily-Borel compactification, Math. Ann. 292, pp. 197 - 240
(1992)

\bigskip\noindent
[R] Ribet K.A., Galois representations attached to eigenforms with nebentypus, In:
Modular functions of one variable V, Springer Lecture Notes in Mathematics 601, pp. 17 - 52 (1977) 

\bigskip\noindent
[JS] Jiang D.-Soudry D., The multiplicity-one theorem for generic automorphic forms
of $GSp(4)$, Pac. J. Math. 229, no. 2, 381-388 (2007)

\bigskip\noindent
[S] Soudry D., The CAP representations of GSp(4),  Crelles Journal 383, pp. 97 -108 (1988)

\bigskip\noindent
[PS] Piateski-Shapiro I.I., On the Saito-Kurokawa lifting, Invent. Math. 71, pp. 309 - 338
(1983)

\bigskip\noindent
[W1] Weissauer R., Four dimensional Galois representations, in: Formes Automorphes (II), Le cas du groupe GSp(4), Editors: Tilouine J., Carayol H., Harris M., Vigneras M.F., Asterisque 302, pp. 67 - 149 (2005)

\bigskip\noindent
[W2] Weissauer R., Endoscopy for GSp(4) and the cohomology of Siegel modular threefolds,
Springer Lecture Notes in Mathematics 1968 (2009)

\bigskip\noindent
[W3] Weissauer R., Existence of Whittaker models  related to four dimensional
symplectic Galois representations, In: Modular Forms on Schiermonnikoog, Editors: 
Edixhofen B., van der Geer G., Moonen B., 
Cambridge University Press pp. 285 - 310 (2008) 

\bigskip\noindent
[T] Taylor R., On the $l$-adic cohomology of Siegel threefolds, Invent. math. 114, pp 289 - 310 (1993)

\bigskip\noindent
[Z] Zagier D., Sur la Conjecture de Saito-Kurokawa, (d'apres Maass), Seminaire de
Theorie de Nombres, Birkh\"auser, Progress in Mathematics, vol. 12, pp. 371 - 394 (1981)

\end{document}